\documentclass[12pt]{amsart}
\usepackage[                
   pdftex,                  
   pdfstartview=FitH,      
   bookmarks=true,          
   bookmarksdepth=section, 
   bookmarksopen=false,     
   bookmarksnumbered=true, 
   ]{hyperref}

\DeclareRobustCommand{\qbinom}{\genfrac[]{0pt}{}}

\newcommand\be            {\begin{equation}}
\newcommand\ee            {\end{equation}}

\theoremstyle{plain}

\newtheorem{theorem}{Theorem}

\newtheorem{corollary}[theorem]{Corollary}
\newtheorem{conjecture}[theorem]{Conjecture}

\theoremstyle{definition}

\numberwithin{equation}{section}
\numberwithin{theorem}{section}
	
\newcommand*{\refh}[2]{\hyperref[#2]{#1~\ref{#2}}} 

\newcounter{ourcount}
\advance\textheight\topskip
\usepackage{amsmath,amsfonts,amsthm,amssymb,amscd}
\usepackage{mathtools} 
\allowdisplaybreaks

 \usepackage{graphicx}
 \usepackage[curve,matrix,arrow,color]{xy}

 \usepackage[usenames,dvipsnames]{color}
 \xyoption{line}

\usepackage[mathscr]{eucal}
\usepackage[OT2,OT1]{fontenc}

\newcommand{\C}{\ensuremath{\mathbb{C}} }

\newcommand{\SL}{\ensuremath{{\mathfrak{sl}}_2}}

\newcommand{\Z}{\mathbb{Z}}

\usepackage{tikz}

\catcode`@=11
\catcode`:=11

\newcount\twistKnot:counter
\def\twistKnot:gobblesign#1{\ifx-#1\else#1\fi}

\def\twistKnot:slopes#1<#2,#3>{\begingroup
	\twistKnot:counter=#1
	\ifnum\twistKnot:counter<0
		\twistKnot:counter=-\twistKnot:counter
		\pgfsetyvec{\pgfpointxy{0}{-1}}%
	\fi
	\ifnum\twistKnot:counter > 0
		\pgfpathcurvebetweentime{#2}{#3}%
			{\pgfpointadd{\pgfpoint\x\y}{\pgfpointxy{-0.1}{0.4}}}%
			{\pgfpointadd{\pgfpoint\x\y}{\pgfpointxy{0.2}{0.4}}}%
			{\pgfpointadd{\pgfpoint\x\y}{\pgfpointxy{0.2}{-0.4}}}%
			{\pgfpointadd{\pgfpoint\x\y}{\pgfpointxy{0.5}{-0.4}}}%
		\pgfpointadd{\pgfpoint\x\y}{\pgfpointxy{0.5}{0}}%
		\pgfgetlastxy\x\y
		\ifnum\twistKnot:counter > 1
			\loop\ifnum\twistKnot:counter > 2
				\pgfpathcurvebetweentime{#2}{#3}%
					{\pgfpointadd{\pgfpoint\x\y}{\pgfpointxy{0}{0.4}}}%
					{\pgfpointadd{\pgfpoint\x\y}{\pgfpointxy{0.2}{0.4}}}%
					{\pgfpointadd{\pgfpoint\x\y}{\pgfpointxy{0.3}{-0.4}}}%
					{\pgfpointadd{\pgfpoint\x\y}{\pgfpointxy{0.5}{-0.4}}}%
				\advance\twistKnot:counter\m@ne\relax
				\pgfpointadd{\pgfpoint\x\y}{\pgfpointxy{0.5}{0}}%
				\pgfgetlastxy\x\y
			\repeat
			\pgfpathcurvebetweentime{#2}{#3}%
				{\pgfpointadd{\pgfpoint\x\y}{\pgfpointxy{0}{0.4}}}%
				{\pgfpointadd{\pgfpoint\x\y}{\pgfpointxy{0.3}{0.4}}}%
				{\pgfpointadd{\pgfpoint\x\y}{\pgfpointxy{0.3}{-0.4}}}%
				{\pgfpointadd{\pgfpoint\x\y}{\pgfpointxy{0.6}{-0.4}}}%
		\fi
	\fi
\endgroup}

\def\twistKnot:twist#1{%
	\ifnum#1>0
		\twistKnot:slopes{#1}<0,1>
			\pgfsetlinewidth{1pt}
			\pgfsetstrokecolor{black}
			\pgfusepath{stroke}
		\twistKnot:slopes{-#1}<0.25,0.75>
			\pgfsetlinewidth{5pt}
			\pgfsetstrokecolor{white}
			\pgfusepath{stroke}
		\twistKnot:slopes{-#1}<0,1>
			\pgfsetlinewidth{1pt}
			\pgfsetstrokecolor{black}
			\pgfusepath{stroke}
	\else
		\twistKnot:slopes{#1}<0,1>
			\pgfsetlinewidth{1pt}
			\pgfsetstrokecolor{black}
			\pgfusepath{stroke}
		\twistKnot:slopes{\twistKnot:gobblesign#1}<0.25,0.75>
			\pgfsetlinewidth{5pt}
			\pgfsetstrokecolor{white}
			\pgfusepath{stroke}
		\twistKnot:slopes{\twistKnot:gobblesign#1}<0,1>
			\pgfsetlinewidth{1pt}
			\pgfsetstrokecolor{black}
			\pgfusepath{stroke}
	\fi
}

\def\twistKnot:connectors#1#2{%
	\pgfpathmoveto
		{\pgfpointadd{\pgfpoint\vx\vy}{\pgfpointxy{0.1}{-0.4}}}%
	\pgfpathcurveto
		{\pgfpointadd{\pgfpoint\vx\vy}{\pgfpointxy{1}{-0.4}}}%
		{\pgfpointadd{\pgfpoint\vx\vy}{\pgfpointxy{1.4}{0.2}}}%
		{\pgfpointadd{\pgfpoint\vx\vy}{\pgfpointxy{1.4}{1.2}}}%
	\pgfpathlineto
		{\pgfpointadd{\pgfpoint\vx\hy}{\pgfpointxy{1.4}{0.6}}}%
	\pgfpathcurveto
		{\pgfpointadd{\pgfpoint\vx\hy}{\pgfpointxy{1.4}{2.6}}}%
		{\pgfpointadd{\pgfpoint\hx\hy}{\pgfpointxy{0.4}{3}}}%
		{\pgfpointadd{\pgfpoint\hx\hy}{\pgfpointxy{0.4}{2.1}}}%
	\pgfpathmoveto
		{\pgfpointadd{\pgfpoint\vx\vy}{\pgfpointxy{0.1}{0.4}}}%
	\pgfpathcurveto
		{\pgfpointadd{\pgfpoint\vx\vy}{\pgfpointxy{0.3}{0.4}}}%
		{\pgfpointadd{\pgfpoint\vx\vy}{\pgfpointxy{0.5}{0.5}}}%
		{\pgfpointadd{\pgfpoint\vx\vy}{\pgfpointxy{0.5}{0.8}}}%
	\pgfpathcurveto
		{\pgfpointadd{\pgfpoint\vx\vy}{\pgfpointxy{0.5}{1.4}}}%
		{\pgfpointxy{0.4}{1}}%
		{\pgfpointxy{0.4}{1.9}}%
}

\newcommand*\doubleTwist[3][]{\begin{tikzpicture}[#1]
	\pgfpointscale{0.25}{\pgfpointxy{\twistKnot:gobblesign#3}{0}}%
	\pgfgetlastxy\vx\vy
	\pgfpointscale{0.5}{\pgfpointxy{0}{\twistKnot:gobblesign#2}}%
	\pgfgetlastxy\hx\hy
	{\pgfpointxy{0}{1}\pgfgetlastxy\x\y
	 \pgfsetyvec{\pgfpointxy{-1}{0}}%
	 \pgfsetxvec{\pgfpoint\x\y}%
	 \pgfpointxy{2}{0}\pgfgetlastxy\x\y
	 \twistKnot:twist{#2}}
	\pgfpoint{-\vx}{\vy}\pgfgetlastxy\x\y
	\twistKnot:twist{#3}%
	\twistKnot:connectors{#2}{#3}%
	\pgfsetxvec{\pgfpointxy{-1}{0}}
	\edef\vx{-\vx}%
	\edef\hx{-\hx}%
	\twistKnot:connectors{#2}{#3}%
		\pgfsetlinewidth{1pt}
		\pgfsetstrokecolor{black}
		\pgfusepath{stroke}
\end{tikzpicture}}

\newcommand*\doubleTwistBoxes[3][]{\begin{tikzpicture}[#1]
	\draw[line width=1pt, solid]
		(0,0) .. controls ++(2, 0) and ++(0,-1.8) ..
		(2,2) .. controls ++(0, 1.8) and ++(0.4, 0) ..
		(1.1,4) .. controls ++(-0.5,0) and ++(0,0.35) .. (0.3,3.4) --
		(0.3,2.2) .. controls ++(0,-1) and ++(0,0.6) ..
		(1.1, 1) .. controls ++(0,-0.3) and ++(0.3,0) .. (0.6,0.6) --
		(-0.6,0.6) .. controls ++(-0.3,0) and ++(0,-0.3) ..
		(-1.1,1) .. controls ++(0,0.6) and ++(0,-1) .. (-0.3,2.2) --
		(-0.3,3.4) .. controls ++(0,0.35) and ++(0.5,0) ..
		(-1.1,4) .. controls ++(-0.4,0) and ++(0,1.8) ..
		(-2,2) .. controls ++(0,-1.8) and ++(-2,0) .. cycle;
	\draw[solid, line width=0.8pt, fill=white, rounded corners=3]
		(-0.6,-0.2) rectangle (0.6,0.8) node[midway] {$#3$}
		(-0.5, 2.2) rectangle (0.5,3.4) node[midway] {$#2$};
\end{tikzpicture}}

\catcode`@=12
\catcode`:=12

\begin{document}

\title{Non-semisimple invariants and Habiro's series}
\author{Anna Beliakova}
\address{University of Zurich, I-Math, Winterthurerstrasse 180, 8008 Zurich}
\email{anna@math.uzh.ch}
\author{Kazuhiro Hikami}
\address{Faculty of Mathematics,
  Kyushu University,
  Fukuoka 819-0395, Japan.}

\email{
  \texttt{khikami@gmail.com}
}
\date{\today}
\begin{abstract} In this paper we establish an explicit  relationship
between Habiro's cyclotomic expansion of the colored Jones polynomial (evaluated at a $p$th root of unity) and the Akutsu-Deguchi-Ohtsuki (ADO) invariants of the double twist knots.  
This allows us to compare the Witten-Reshetikhin-Turaev (WRT) and Costantino-Geer-Patureau (CGP) invariants
of 3-manifolds obtained by $0$-surgery on these knots.
The difference between them is determined by the $p-1$
coefficient of the Habiro series. We expect these to hold for all  Seifert genus 1 knots.
\end{abstract}





\maketitle

\section{Introduction}
In \cite{Ha}
 Habiro stated the following result.
 Given a 0-framed knot $K$ and an $N$-dimensional irreducible representation of the quantum $\SL$,  there  exist polynomials
 $C_n(K;q)\in \Z[q^{\pm 1}]$, $n\in {\mathbb{N}}$, such that
\begin{equation}\label{Jones}
  J_K(q^N,q) = \sum_{n=0}^\infty C_n(K;q) \,
  (q^{1+N};q)_n ( q^{1-N}; q)_n
\end{equation}
is the $N$-colored Jones polynomial of $K$.
Here $(a;q)_n=(1-a)(1-aq)\dots(1-aq^{n-1})$ and $q$ is a generic parameter.  Replacing $q^N$ by a formal variable $x$ we get 
\begin{align}
\label{cyclotomic}
  J_K(x,q) = \sum_{n=0}^\infty C_n(K;q) \,
  (xq;q)_n ( x^{-1}q; q)_n =
  \sum_{m\geq 0} a_m(K;q)\, \sigma_m  (x,q)
\end{align}
where
$$\sigma_m(x,q)=\prod_{i=1}^m\left(x+x^{-1} -q^i-q^{-i} \right) \quad \text{and}\quad
C_m(K;q)= (-1)^m q^{-\frac{1}{2} m(m+1)} a_m(K;q) 
 $$
 known as {\it cyclotomic expansion}
of the colored Jones polynomial or simply  {\it Habiro's series}. This expression dominates all colored Jones polynomials and converges
in the cyclotomic completion
$$ {\lim\limits_{\overleftarrow{\hspace{2mm}n\hspace{2mm}}}}\;   
 \Z[q^{\pm 1}][x+x^{-1}]/(\sigma_n(x,q))$$
 of the  center of  $U_q(\SL)$. In details, \eqref{cyclotomic} belongs to the cyclotomic
completion of the even part of the center after the
 identification 
 $C^2-2=x+x^{-1}$, where $C$ is the Casimir.
 
Habiro's series played a central role in the construction of the unified invariants of homology 3-spheres
\cite{Ha, HL, BBL} dominating all WRT invariants.
In \cite{BCL} it was used to prove  integrality of 
the WRT invariants for all 3-manifolds at all roots of unity.

Given the power of Habiro's series, it is not surprising that they are notoriously difficult to compute.
 So far,  Habiro's cyclotomic expansions were computed  explicitly 
for 
the  following   infinite families:
   the twist knots~\cite{Ma}, the $(2,2t+1)$ torus knots~\cite{KL},
   and the double twist knots~\cite{LO}.

Recently the non-semisimple quantum invariants of links and 3-manifolds attracted a
 lot of attention. Physicists expect them to play a crucial role in categorification of quantum 3-manifold invariants \cite{Gukov}. Mathematicians resolved the problem of
 nullity of these invariants by introducing the
 {\it modified traces} \cite{GPT,BBGa}.

 The aim of this paper is to connect the Habiro cyclotomic expansion with  the non-semisimple world.
The non-semisimple  invariants arise 
in 
  specializations of the quantum $\SL$ at $q=e_p$, the primitive $p^{\text{th}}$  root of unity. 
 The
ADO link invariant is
 obtained in the setting of the unrolled
quantum $\mathfrak{sl}_2$  \cite{unrolled}.
This group admits   $p$-dimensional irreducible projective modules
$V_\lambda$  whose highest weights $\lambda$ are given by any 
complex number.
Even through the definition of the unrolled quantum group
requires a choice of the 
square  root of  $e_{p}$ (that we denote by
$e_{2p}$), the ADO invariant of a 0-framed knot
does not depend on this choice.
The representation category of the unrolled $\SL$ is ribbon, and hence, the ADO invariant
can be defined by applying the usual Reshetikhin-Turaev construction, i.e.
we
color the $(1,1)$-tangle $T$ whose closure is $K$
with  $V_\lambda$, the Reshetikhin-Turaev functor sends then $T$ to an endomorphism
of $V_\lambda$  that is ${\mathrm{ADO}_K}(e_p^{\lambda+1}, e_{p}) \text{id}_{V_\lambda}$.
   
  The corresponding non-semisimple 3-manifold invariant 
for a pair $(M,\lambda)$, where $M$ is a closed oriented 3-manifold and $\lambda\in H^1(M, \C/2\Z)$ is a cohomology class, was defined
  in \cite{CGP}. 
If $M=S^3(K)$ is obtained by surgery on a knot $K$ in $S^3$
with non-zero framing, then
$M$ is a rational homology 3-sphere and $\lambda$ is rational. In this case, the
 CGP invariant of $M$ was shown  to be determined by the Witten-Reshetikhin-Turaev invariant (WRT) in \cite{WRTCGP}.   It remains to analyze the case of $0$-framed surgeries with $\lambda \neq 0,1$.
For a $0$-framed knot $K$,
\begin{equation*}
  \mathrm{CGP}(S^3(K), \lambda)=
  \sum^{p-1}_{n=0}
  d^2(\lambda+2n) \mathrm{ADO}_K (e_p^{\lambda+2n+1}, e_{p})
\end{equation*}
where $ d(\lambda+2n)$ is the modified dimension of $V_{\lambda+2n}$.
   
\section{Our results}
Recently, in \cite{Will} Willetts constructed the knot invariant $$F_\infty(q, x;K)\in \;\hat R:=\;
{\lim\limits_{\overleftarrow{\hspace{2mm}n\hspace{2mm}}}}\;   \;
\frac{ \Z[q^{\pm 1/2}, x^{\pm 1/2}]}{I_n}$$
where $I_n$ is the ideal generated by $$\left\{ \,\prod^n_{i=0}\,
(x^{\frac{1}{2}} q^{\frac{l+i}{2}}-x^{-{\frac{1}{2}}}q^{-\frac{l+i}{2}})\,|\; l\in \Z\right\}$$
  that dominates the colored Jones polynomials and the ADO
invariants of $K$.
In details,
\begin{equation}
  \label{ADO_F}
 F_\infty(q,q^N;K)=J_K( q^N,q)\quad\text{and}\quad
 F_\infty(e_p, x;K)=\frac{\mathrm{ADO}_K(x, e_{p})}{\Delta_K(x^p)}
\end{equation}
where $\Delta_K(x)$ is the Alexander polynomial of  $K$. 
The famous
Melvin-Morton-Rozansky theorem \cite{Ro} follows from the above result at $p=1$.

We claim that $F_\infty(q,x;K)$ coincide with the
Habiro expansion~\cite{Ha,HL} $J_K(x,q)$ as elements of $\hat R$. 
Hence, the result of Willetts in \cite{Will} can be reformulated as follows.

\begin{theorem}\label{main1}
The universal Habiro series determine the ADO invariants, i.e.
 $$J_K( x, e_{p})= \frac{{\mathrm{ADO}_K}(x, e_{p})}{\Delta_K(x^p)}
.$$
\end{theorem}



\vspace{1mm}

\noindent
{\bf Example.} For the first two knots at the first two roots of unity this looks as follows:
$$
J_{3_1}(t,q)=\sum_{m {\geq 0}} (-1)^m q^{m(m+3)/2} \,\sigma_m(t,q),
\quad
J_{3_1}(t,1)=\sum_{m {\geq 0}} (-1)^m(t+t^{-1} -2)^m=\frac{1}{t+t^{-1}-1},$$
$$J_{3_1}(t,-1)=(1-(t+t^{-1}+2))\sum_{m{\geq 0}} (-1)^m (t^2+t^{-2} -2)^m=\frac{\Delta_{3_1}(-t)}{\Delta_{3_1}(t^2)},$$
$$J_{4_1}(t,q)=\sum_{m{\geq 0}} \sigma_m(t,q),
\quad J_{4_1}(t,1) =\sum_{m {\geq 0}} (t+t^{-1} -2)^m=\frac{1}{1-(t+t^{-1}-2)}=\frac{1}{\Delta_{4_1}(t)},
$$
$$J_{4_1}(t,-1)=(1+t+t^{-1}+2)\sum_{m {\geq 0}}  (t^2+t^{-2} -2)^m=\frac{\Delta_{4_1}(-t)}{\Delta_{4_1}(t^2)}.$$
\vspace{1mm}

A challenging open problem is to find an explicit formula
for ADO using the coefficients $a_n(K;e_p)$ of the  Habiro series.
In the examples above this is possible due to  periodicity of these coefficients.
Our next 
theorem  
generalizes these examples to all double twist knots.

Let  $\{K_{(l,m)}\,|\, l,m\in \Z\}$ be the 2-parameter family of double
twist
knots such that 
$K_{(l,m)}=K_{(m,l)}$, $K_{(1,1)}=3_1$ and
$K_{(-1,1)}=4_1$   depicted in Figure 1. 

\begin{figure}
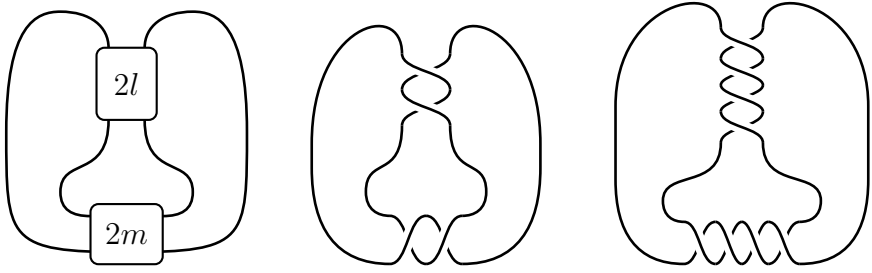

\doubleTwistBoxes[x=8mm,y=8mm]{2l}{2m}\quad\quad
\doubleTwist[x=8mm,y=8mm]{2}{2} \quad\quad
\doubleTwist[x=7mm,y=7mm]{4}{-4}
\caption{The double twist knots: $K_{(l,m)}$, the trefoil and $K_{(2,-2)}$. The integers in the boxes denote the number of  half twists.}
\end{figure}

\begin{theorem}\label{main2}
Let $K=K_{(l,m)}$ be a  double twist knot. Then we have a new 
expression for its
ADO invariant
$$\mathrm{ADO}_{K}(x, e_{p})=\sum^{p-1}_{n=0}
a_n(K;e_p) \,\sigma_n( x, e_p)\; .
$$
 In addition, for all natural numbers $n,k$ we have
$$a_{n+kp}(K;e_p)=a_n(K;e_p)\; a_k(K;1)\quad
\text{and} \quad\Delta^{-1}_K(t)=\sum_k a_{kp}(K;e_p)(t+t^{-1}-2)^k .
$$
\end{theorem}

Our next result establishes a relationship between
CGP and WRT invariants of 3-manifolds obtained by $0$-surgery on double twist knots.  Note that for rational surgery
on knots CGP is known to be determined by WRT and the order
of the first homology group \cite{WRTCGP}.

\begin{theorem}\label{main3}
Let $M=S^3(K_{(l,m)})$ where $K_{(l,m)}$ is 0-framed double twist knot and $\lambda\in \C/2\Z$ with $\lambda \neq 0,1$.
 For odd $p>1$ we have
$$\mathrm{CGP}(M,\lambda)= \frac{ 1}{{\{p\lambda\}}^2} \mathrm{WRT}(M) +p\, a_{p-1}(K_{(l,m)};e_p)
\quad\text{
where} \quad \{y\}=e_{2p}^{y}-e_{2p}^{-y}.$$
\end{theorem}

 \begin{corollary} The coefficient
$a_{p-1}(K_{(l,m)};e_p)$ of the Habiro series 
is a topological invariant of the 3-manifold 
obtained by $0$-surgery on the double twist knot $K_{(l,m)}$.
\end{corollary}

For $p=2$,  $\mathrm{WRT}(M)=1$ and $\mathrm{CGP}(M)$ 
 is the Reidemeister torsion of $M$ according to
\cite[Thm. 6.23]{BCGP}.
Hence, one can think about the coefficient $a_{p-1}(K;e_p)$
as a generalization of the Reidemeister torsion.
It is an interesting open problem to find a topological interpretation of this invariant. In general we would expect following to hold.

\begin{conjecture} Theorems  \ref{main2}, \ref{main3} hold for
any Seifert genus 1 knot.
\end{conjecture}

Examples of knots with higher Seifert genus
provide $(2,2t+1)$ torus knots
with $t\geq 2$. Our next result is a computation of
the ADO invariants for this family of knots.
\begin{theorem}
  \label{ado_2t+1}
  Let $K=T_{(2,2t+1)}$ be a torus knot.
  Then the ADO invariant is
 \begin{equation*}
  \mathrm{ADO}_{K}(x,e_p)
  =
  e_p^{~ t} x^{t(1-p)}
  \sum_{k_t\geq \dots \geq k_1 \geq 0}^{p-1}
  (x \, e_p ; e_p)_{k_t} \,  x^{k_t}
  \prod_{i=1}^{t-1} e_p^{~ k_i(k_i+1)} x^{2k_i}
  \begin{bmatrix}
    k_{i+1} \\
    k_i
  \end{bmatrix}_{e_p}
\end{equation*}
\end{theorem}
An interesting open problem is to determine periodicity
of $\{a_n(T_{(2,2t+1)},e_p)\}_{n\geq 0}$ analogous to those for double twist knots in 
Theorem \ref{main2}.
In Appendix we show that $\{a_n(T_{(2,5)},e_p)\}_{n\geq 0}$ do not satisfy  properties listed in 
Theorem \ref{main2}.

In the last section we compute CGP and WRT invariants for $0$-surgeries on  $(2,2t+1)$ torus knots and compare them.
We observe that  (up to normalization) CGP can be viewed
as a Laurent polynomial in $T^{\pm 1}:=(e^{\pm \lambda}_p)^p$, such that
its evaluation at $T=1$  reproduces  the WRT invariant. 
In contrast to Theorem \ref{main3} however,
 for torus knots the WRT does not form anymore the degree zero part of this Laurent polynomial.

\bigskip
\noindent
\textbf{Acknowledgement} AB would like to thank Christian Blanchet for many helpful discussions, and Krzysztof Putyra                                    
for providing  pictures of the double twist knots.
The work of KH is partially supported by JSPS KAKENHI Grant Numbers,
JP16H03927, JP20K03601, JP20K03931.

\section{Proofs} In this section we prove our four Theorems.

\subsection{Preliminaries}
We set
$[n]_q=\frac{1-q^n}{1-q}$,
$[n]_q!=[n]_q \dots [2]_q [1]_q$.
The $q$-binomial is defined by
$$ \begin{bmatrix}
  n \\ k
\end{bmatrix}_q=\frac{[n]_q!}{[k]_q![n-k]_q!}
$$
Observe that evaluating at $q=e_p$, the primitive $p^{\text{th}}$ root
of unity,
we have
 \begin{equation}\label{fact:qbinom}
  \begin{bmatrix}
    n + a p \\ k +bp
  \end{bmatrix}_{e_p}=\begin{bmatrix}
    n   \\ k 
  \end{bmatrix}_{e_p} {a\choose b}
  \end{equation} 
where the last factor is the usual binomial coefficient. 
Moreover, setting $y=1$ in  the identity
$$x^p+x^{-p}-y^p-y^{-p}=\prod^{p}_{i=1} \left(x+x^{-1}-
ye^i_p -y^{-1}e^{-i}_p\right),$$
we get the following   relation 
\begin{equation}\label{sigmap}
 \sigma_p(x, e_p)=x^p+x^{-p}-2=(1-x^p)(x^{-p}-1)
 \end{equation}

\subsection{Proof of Theorem \ref{main1}}
For $J_K(x,q)$ in~\eqref{cyclotomic}, it is easy to check that its evaluations  
at $x=1,q, q^2, \dots$ determine the coefficients
$C_n(K;q)$ recursively.
Explicitly,
\begin{equation*}
  C_n(K;q)
  =
  -q^{n+1}\sum_{l=1}^{n+1}
  \frac{
    (1-q^l) (1-q^{2l})}{
    (q)_{n+1-l}(q)_{n+1+l}
  }
  (-1)^l q^{\frac{1}{2}l(l-3)}
  J_K(q^{l},q)\ .
\end{equation*}
Moreover, 
$J_K(x,q) \in \hat R$ by \eqref{cyclotomic}.
In   \cite[Prop. 57]{Will} Willetts proved
that the
 evaluations at $x=q^N$ for $N\in \mathbb N$ determine 
   $F_\infty(q,x;K)$.  Moreover,  both invariants
   coincide
  at all evaluations, since  $ F_\infty(q,q^N;K)=J_K( q^N,q)$. Hence,
  $$F_\infty(q,x;K)=J_K(x,q) \in \hat R$$ 
and the result follows from the Willetts theorem.
$\hfill\Box$

\subsection{Proof of Theorem \ref{main2}}



Let $K=K_{(l,m)}$ be a double twist knot.
For $p=1$,  Theorem \ref{main1} reads $J_K(t,1)=\frac{1}{\Delta_K(t)}$
and hence
{using~\eqref{sigmap}}
$$\frac{1}{\Delta_K(x^p)}=\sum^\infty_{k=0} a_k(1) \sigma^k_p$$
where we set $\sigma_k:=\sigma_k(x, e_p)$ and $a_k(e_p):=a_k(K;e_p)$
for brevity.
From Theorem \ref{main1} assuming \eqref{ADOtwist}
we get
\be\label{fact}
J_K(x,e_p)=\sum^{\infty}_{k=0}\left(
\sum^{p-1}_{n=0} a_{n+pk}(e_p)\; \sigma_n
\right) \sigma^k_p= 
\sum^{p-1}_{n=0} a_n(e_p) \sigma_n\;
\sum^\infty_{k=0} a_k(1)\sigma^k_p  .
\ee
Here we used 
$\sigma_{n+kp}=\sigma_n
\sigma^k_p$.
The desired symmetries for the Habiro coefficients follow.
It remains to
 prove
\begin{equation}\label{ADOtwist}
\mathrm{ADO}_{K}(e^\lambda_p, e_{p})=\sum^{p-1}_{n=0}
a_n(e_p) \,\sigma_n(  e_{p}^\lambda, e_p)\; .
\end{equation}


For this purpose, let us first observe that
the Alexander polynomial of the double twist knot $K$ is
(see, \emph{e.g.},~\cite{Hill})
\begin{equation}
  \label{eq:10}
  \begin{aligned}
    \Delta_{K}(x)
    & =
    1+lm(x+x^{-1}-2)
    \\
    & = f_{lm}( - \tfrac{(1-x)^2}{x})
  \end{aligned}
\end{equation}
where
  $f_n(z)=1-n z$. Hence, 
\begin{equation*}
  \frac{1}{f_n(z)}=1+n z + n^2 z^2 +n^3 z^3 + \dots
\end{equation*}
 and 
\begin{equation}
  \label{eq:11}
  \frac{1}{\Delta_{K}(x)} = 
  \sum_{n=0}^\infty {(-lm)^n}\;(x+x^{-1}-2)^{n}=
  \sum_{n=0}^\infty {l^n m^n}\;(1-x)^n(1-x^{-1})^{n} 
\end{equation}

The Habiro expansion of the colored Jones polynomial for double twist knots is given by Lovejoy and Osburn \cite{LO}. For $l,m>0$ we have
\begin{gather}
  \label{U_twist_plus}
  J_{K_{(l,m)}}(x,q)
  =
  \sum_{\substack{n=t_m \geq \dots\geq t_1 \geq 0\\
      n=s_l \geq \dots \geq s_1 \geq 0}}
  (x q;q)_{n}(x^{-1}q;q)_{n}
  \,
  q^{n}
  \,
  \prod_{i=1}^{l-1}
  q^{s_i(s_i+1)}
  \begin{bmatrix}
    s_{i+1}\\ s_i
  \end{bmatrix}_{q}
  \prod_{j=1}^{m-1}
  q^{t_j(t_j+1)}
  \begin{bmatrix}
    t_{j+1}\\ t_j
  \end{bmatrix}_{q}
  \\
  \label{U_twist_minus}
  J_{K_{(l,-m)}}(x,q)
  =
  \sum_{\substack{
      n=t_m \geq \dots\geq t_1 \geq 0
      \\
      n=s_l \geq \dots \geq s_1 \geq 0}}
  (x q;q)_{n}(x^{-1}q;q)_{n}
  \,
  (-1)^{n} q^{-\frac{1}{2} n(n+1)}
  \prod_{i=1}^{l-1}
  q^{s_i(s_i+1)}
  \begin{bmatrix}
    s_{i+1} \\ s_i
  \end{bmatrix}_{q}
  \prod_{j=1}^{m-1}
  q^{-t_j(t_{j+1}+1)}
  \begin{bmatrix}
    t_{j+1}\\ t_j
  \end{bmatrix}_{q}
\end{gather}
Note that due to the following symmetries: 
$K_{(l,m)}=K_{(m,l)}$ and $K_{(-l,-m)}$ is the mirror image of $K_{(l,m)}$, the expressions above cover all double twist knots (up to substitution of $q$ by $q^{-1}$).

Evaluating \eqref{U_twist_plus} at $q=e_p$ we get
\begin{equation*}
  \begin{aligned}
    &J_{K_{(l,m)}}(x, e_p )
    \\
    = &
    \sum_{k=0}^\infty
    \sum_{n=0}^{p-1}(xe_p; e_p)_{k p +n}
    (x^{-1}e_p; e_p)_{k p +n}  e_p^{ n}
    \\
    & \qquad \times
    \sum_{\substack{
    k p +n =t_m \geq t_{m-1}
      \geq \dots \geq
      t_1 \geq 0
      \\
      k p +n =s_l \geq s_{l-1}
      \geq \dots \geq
      s_1 \geq 0}}
    \prod_{i=1}^{l-1}
    e_p^{~ s_i(s_i+1)}
    \begin{bmatrix}
      s_{i+1} \\
      s_i
    \end{bmatrix}_{e_p}
    \prod_{j=1}^{m-1}
    e_p^{~ t_j(t_j+1)}
    \begin{bmatrix}
      t_{j+1} \\
      t_j
    \end{bmatrix}_{e_p}
    \\
    =&
    \sum_{k=0}^\infty (1-x^p)^k (1-x^{-p})^k \, (l m)^k\, \sum_{n=0}^{p-1}
    (xe_p ; e_p)_n
    (x^{-1}e_p; e_p)_n \,
    e_p^{ n}  \\
    &
    \qquad \times
    \sum_{\substack{
    n =t_m \geq t_{m-1}
      \geq \dots \geq
      t_1 \geq 0\\
      n =s_l \geq s_{l-1}
      \geq \dots \geq
      s_1 \geq 0}}
    \prod_{i=1}^{l-1}
    e_p^{~ s_i(s_i+1)}
    \begin{bmatrix}
      s_{i+1} \\
      s_i
    \end{bmatrix}_{e_p}
    \prod_{j=1}^{m-1}
    e_p^{~ t_j(t_j+1)}
    \begin{bmatrix}
      t_{j+1} \\
      t_j
    \end{bmatrix}_{e_p}
    \\
    =&
    \frac{1}{\Delta_{K_m}(x^p)}
    \sum_{n=0}^{p-1} (x e_p; e_p)_n
    (x^{-1}e_p; e_p)_n \,
    e_p^{n}
    \sum_{
      \substack{
        n =t_m \geq t_{m-1}
        \geq \dots \geq
        t_1 \geq 0\\
        n =s_l \geq s_{l-1}
      \geq \dots \geq
      s_1 \geq 0}}
\prod^{l-1}_{i=1}   
    \begin{bmatrix}
      s_{i+1} \\
      s_i
    \end{bmatrix}_{e_p}
    \prod_{j=1}^{m-1}
    e_p^{~ t_j(t_j+1)}
    \begin{bmatrix}
      t_{j+1} \\
      t_j
    \end{bmatrix}_{e_p}
  \end{aligned}
\end{equation*}
Here we have used
\eqref{fact:qbinom} and \eqref{sigmap} (see also \eqref{Andrews_root_unity}).
Thus applying Theorem \ref{main1} we get \eqref{ADOtwist}.
%
The case of $K_{(l,-m)}$ is similar.
As a result, for any double twist knot $K$
we have
$$ \mathrm{ADO}_{K}( x, e_p)= \sum^{p-1}_{n=0} a_n(K;e_p) \sigma_n(x,e_p).$$
$\hfill\Box$

\subsection{Proof of Theorem \ref{main3}}

The non-semisimple invariant is defined as follows
\begin{equation}
  \label{CGP_ADO}
  \mathrm{CGP}(S^3(K), \lambda)=
  \sum^{p-1}_{n=0}
  d^2(\lambda+2n) \mathrm{ADO}_K (e_p^{\lambda+2n+1}, e_{p})
\quad\text{
  where} \quad d(y)=\frac{\{y+1\}}{\{py\}}
\end{equation}
is the modified dimension
and $\{y\}=e^y_{2p} -e^{-y}_{2p}$.
Inserting the new expression for the ADO invariant and
exchanging the sums we get
\begin{equation}\label{CGP}
\mathrm{CGP}(S^3(K), \lambda)= \frac{1}{\{p\lambda\}^2}
\sum^{p-1}_{m=0} a_m(e_p) \sum^{p-1}_{n=0} \{\lambda+2n+1\}^2
\sigma_m(e^{\lambda+2n+1}_{p}, e_p)
\end{equation}
On the other hand, for odd $p$ the WRT invariant can be written as
follows
{
  (see~\cite{Le})
}
\begin{equation}\label{WRT}
  \begin{aligned}[t]
    \mathrm{WRT}(S^3(K))
    & =
    \sum_{\substack{
        0<n<2p\\
        \text{$n$: odd}
      }}
    \left( e_{2p}^n - e_{2p}^{-n} \right)^2
    J_K(x=e_p^{-n}, e_p)
    \\
    &=
    \sum^{(p-3)/2}_{m=0} a_m(e_p) \sum^{p-1}_{n=0} \{2n+1\}^2
    \sigma_m(e^{2n+1}_{p}, e_p)
  \end{aligned}
\end{equation}
(see,\emph{e.g.}~\cite{Le}). 
 The usual normalization of the WRT can be obtained by multiplying with $\{1\}^{-2}$.
Both invariants can be computed explicitly using so-called Laplace transform method. 
 Observe that up to normalization both invariants can be written as
\begin{equation*}
  \sum_{m} a_m(e_p) \left(\sum^{p-1}_{n=0} \{z+m\}
    \{z+m-1\}\dots\{z+1\} \{z\}^2 \{z-1\}\dots\{z-m+1\}
    \{z-m\}\right)
\end{equation*}
where $z=\lambda+2n+1$ for the CGP
and $z=2n+1$ for the WRT invariants. The expression
in the brackets
  is a monic polynomial of degree $m+1$
in $e_p^z+ e_p^{-z}$.  
 Moreover, for 
an odd  $p$ and
any { $a\in \mathbb{Z}$ s.t. $0\leq |a|\leq p$}
we have
\begin{equation}
  \label{sum_ep}
  \sum^{p-1}_{n=0} e_p^{(\lambda +2n+1)a}=
  \begin{cases}
    0 &\text{if $ a\neq 0, \pm p$}\\
    p &\text{if $ a=0$}\\
    pe_p^{\pm  \lambda p}&\text{if $a=\pm p $}
  \end{cases}
\end{equation}
The contribution 
from terms with $m< (p-1)/2$ and $a=0$ in CGP
is exactly  the WRT invariant
of $M=S^3(K)$, i.e.
\begin{equation*}
  \mathrm{WRT}(M)
  =
  -2p\sum_{m=0}^{(p-3)/2}(-1)^m a_m(e_p)
  \begin{bmatrix}
    2m+1 \\
    m
  \end{bmatrix}_{e_p}
  e_p^{-\frac{1}{2}m(m+1)}
\end{equation*}
(compare  \cite[Prop.3.1]{Le}).
The next coefficients for $a=0$ and $(p-1)/2 \leq m<p-1$ are zero, since the $q$-binomial  in
this case contains $\{p\}=0$.
The contribution for $a=\pm p$ and $m=p-1$ gives
$pa_{p-1}(e_p^{\lambda p}+e_p^{-\lambda p})$ and for
$a=0$ and $m=p-1$ 
$$
(-2)p\, a_{p-1}(e_p)\, (-1)^{p-1}
 {\qbinom{2p-1}{p}}_{e_{2p}}=-2pa_{p-1}(e_p).$$
Using that $\{p\lambda\}^2=e_p^{\lambda p}+e_p^{-\lambda p}-2$ we get the result. 
$\hfill\Box$

\subsection{Proof of Theorem~\ref{ado_2t+1}}

  It is known that the $N$-colored Jones polynomial for knot $K$ satisfies
  the recurrence relation.
  For instance, for torus knot $K=T_{(s,t)}$,
  it is~\cite{KH04}
  \begin{multline}
    \label{recurrence_torus}
    J_{T_{(s,t)}}(q^N,q)
    =
    \frac{q^{\frac{1}{2}(s-1)(t-1)(1-N)}}{1-q^{-N}}
    \left(
      1-q^{s(1-N)-1} -q^{t(1-N)-1} + q^{(s+t)(1-N)}
    \right)
    \\
    +\frac{1-q^{2-N}}{1-q^{-N}} q^{st(1-N)-1}
    J_{T_{(s,t)}}(q^{N-2},q)
  \end{multline}

  For a case of $K=T_{(2,2t+1)}$,
  we have a $q$-hypergeometric series~\cite{KH02}
  \begin{equation}
    \label{Jones_2t+1}
    J_{K}(x,q)
    =
    (q x)^t
    \sum_{k_t \geq \dots \geq k_2 \geq k_1\geq 0}
    (q x; q)_{k_t}
    x^{k_t}
    \prod_{i=1}^{t-1}
    q^{k_i(k_i+1)}x^{2k_i}
    \begin{bmatrix}
      k_{i+1} \\
      k_i
    \end{bmatrix}_q \in \hat{R},
  \end{equation}
  where the $N$-colored Jones polynomial is given by setting $x=q^{-N}$.
  Following~\cite{Will}, we get the ADO invariant by putting
  $q=e_p$ in this expression and then by multiplying the result with the Alexander polynomial.
   Note that 
\begin{equation*}
  \Delta_{K}(x)
  =
  x^{-t} \frac{1+x^{2t+1}}{1+x} .
\end{equation*}

  As an application of~\eqref{fact:qbinom},
we have
\begin{equation}
  \label{Andrews_root_unity}
  \sum_{p+m=k_t\geq \dots \geq k_1 \geq 0}
  \prod_{i=1}^{t-1} e_p^{~ k_i(k_i+1)} x^{2k_i}
  \begin{bmatrix}
    k_{i+1} \\
    k_i
  \end{bmatrix}_{e_p}
  =
  \frac{1-x^{2 t p}}{1-x^{2p}}
  \sum_{m=k_t\geq \dots \geq k_1 \geq 0}
  \prod_{i=1}^{t-1} e_p^{~ k_i(k_i+1)} x^{2k_i}
  \begin{bmatrix}
    k_{i+1} \\
    k_i
  \end{bmatrix}_{e_p} .
\end{equation}
Then we see that
\begin{equation*}
  J_{K}(x, e_p)
  =
  \frac{1}{\Delta_{K}(x^p)}
  e_p^{~ t} x^{t(1-p)}
  \sum_{k_t\geq \dots \geq k_1 \geq 0}^{p-1}
  (x e_p ; e_p)_{k_t} x^{k_t}
  \prod_{i=1}^{t-1} e_p^{~ k_i(k_i+1)} x^{2k_i}
  \begin{bmatrix}
    k_{i+1} \\
    k_i
  \end{bmatrix}_{e_p}
\end{equation*}
This implies the statement of the theorem.
$\hfill\Box$

\vspace{2mm}

\noindent
{\bf Example.} In the case of $p=2$,
\emph{i.e.} $e_2=-1$,
we have
$
\left[
  \begin{smallmatrix}
    a \\ b
  \end{smallmatrix}
\right]_{-1}
=1
$
when $0 \leq b \leq a \leq 1$.
Then we get the known result that the ADO polynomial coincides with
the Alexander
polynomial as follows.
\begin{equation*}
  \begin{aligned}[b]
    \mathrm{ADO}_K(x,-1)
    & =
    (-1)^t x^{-t}
    \sum_{k_t \geq \dots \geq k_1 \geq 0}^1
    (-x;-1)_{k_t} x^{k_t} x^{2(k_1+\dots+k_{t-1})}
    \\
    & =
    (-x)^{-t}
    \left(
      1+(1+x) x \frac{1-x^{2t}}{1 - x^2}
    \right)
    \\
    & =
    \Delta_{K}(-x)
  \end{aligned}
\end{equation*}


\vspace{3mm}

\noindent
{\bf Remark on Theorem~\ref{ado_2t+1}.}
Setting  $x=q^{-N}$ in the recurrence relation~\eqref{recurrence_torus}, 
we get
\begin{equation*}
  J_{K}(x,q)
  =
  \frac{(x q)^{\frac{1}{2}(s-1)(t-1)}}{1-x}
  \left(
    1-q^{s-1}x^s -q^{t-1}x^t + q^{s+t} x^{s+t}
  \right)
  +\frac{1-x q^{2}}{1-x} q^{st-1} x^{st}
  J_K(x q^2, q)
\end{equation*}
whose solution is a $q$-series
as
  \begin{equation}
    J_{T_{(s,t)}}(x,q)
    =
    \frac{( q x )^{\frac{1}{2}(s-1)(t-1)}}{1-x}
    \sum_{n\geq 0}
    \chi_{s,t}(n)\,
    q^{\frac{n^2-(s t -s -t)^2}{4 st}}
    x^{\frac{1}{2}(n-(s t - s -t))}
  \end{equation}
  where
  \begin{equation*}
    \chi_{s,t}(n)=
    \begin{cases}
      1, & \text{for $n=s t \pm(s+t) \mod 2 st$},
      \\
      -1, & \text{for $n= s t \pm(s-t) \mod 2 s t$},
      \\
      0, & \text{otherwise.}
    \end{cases}
  \end{equation*}
  It is conjectured~\cite{GHNPPS} that this gives the ADO invariant
  \begin{equation}
    \label{ado_st}
    \begin{aligned}[t]
      \mathrm{ADO}_{T_{(s,t)}}(x, e_p)
      & =
      \Delta_{T_{(s,t)}}(x^p) \, J_{T_{(s,t)}}(x, e_p)
      \\
      & =
      \frac{x^{\frac{1}{2} - \frac{1}{2} (s-1)(t-1) p}
        (1-x^p)}
      {(1-x)(1-x^{s p})(1-x^{t p})}
      e_p^{\frac{1}{4} ( s t - \frac{s}{t} - \frac{t}{s})}
      \sum_{l=0}^{2 s t p}
      \chi_{ s, t }(l)
      e_p^{\frac{l^2}{4 s t}} x^{\frac{l}{2}}
    \end{aligned}
  \end{equation}
  Numerical computations for some $p$'s in the case of $T_{(2,2t+1)}$
  support the equality,
  Theorem~\ref{ado_2t+1} 
  and~\eqref{ado_st}.

    For
  a mirror image of $T_{(2,2t+1)}$,
  the coefficients of the Habiro expansion~\eqref{cyclotomic} were determined in~\cite{KL} as
  \begin{equation}\label{Hab.coef}
    a_n(
    \overline{T_{(2,2t+1)}};q)
    =
    (-1)^n q^{\frac{1}{2}n(n+1)
      +
      n+1-t}
    \sum_{n+1=k_t \geq k_{t-1} \geq \dots \geq k_1 \geq 1}
    \prod_{i=1}^{t-1}
    q^{k_i^2}
    \begin{bmatrix}
      k_{i+1} + k_i-i + 2\sum_{j=1}^{i-1} k_j
      \\
      k_{i+1}-k_i
    \end{bmatrix}_q\, .
  \end{equation}
  This is written in $q$-binomial, and can be evaluated at  $q=e_p$.
  Though, contrary to the double twist knots, 
  $a_n(K;e_p)$ do not have a simple periodicity.
  See Appendix.

\section{CGP versus WRT for $0$-surgeries on torus knots}
Throughout this section $K=T_{(2,2t+1)}$ is a torus knot and $M$ is a 3-manifold  obtained by $0$-surgery on $K$. 
We can compute
the WRT invariant of $M$ by inserting~\eqref{Jones_2t+1}
into the definition
\begin{equation}
  \label{eq:31}
  \mathrm{WRT}
 (M)
  =
  \sum_{\substack{
      0<n<2p\\
      \text{$n$:odd}
    }}
  \left.
    \left(e^n_{2p}-e^{-n}_{2p}
          \right)^2
    J_K(q^{-n}; e_p)
  \right|_{q=e_p}
\end{equation}
which gives
\begin{equation}
  \label{eq:WRTtorus}
  \begin{aligned}[b]
    \mathrm{WRT}(M)
    &=
    \frac{
      e_p^{t}}{2}\;
    \sum_{k=0}^{p-1} (-1)^k e_p^{\frac{2t+1}{2}k^2+\frac{2t-1}{2}k}
    \sum_{\substack{
        0<n<2p\\
        \text{$n$:odd}
      }}
    e_p^{\left( 1 - t-(2t+1)k \right)n}
    \left(1-e_p^{-n} \right)
    \left( 1-e_p^{2k+1} e_p^{-2n}\right)
    \\
    &
    =
    \frac{1}{2}
    \sum_{k=0}^{p-1}
    (-1)^k
    e_p^{\frac{2t+1}{2}k^2+\frac{2t-1}{2}k
      -(2t+1)k
    }
    \sum_{n=0}^{p-1}
    e_p^{-2n \left( t+(2t+1)k \right)}
    \left( e_p^{2n+1}-1 \right)
    \left( 1-e_p^{2k-4n-1} \right)
  \end{aligned}
\end{equation}

The CGP
invariant of $M$ can be computed using 
the ADO given in Theorem~\ref{ado_2t+1}.
A
simpler expression of the CGP invariant can be obtained from~\eqref{ado_st}.
In our case,
~\eqref{ado_st} can be rewritten as
\begin{equation*}
  \mathrm{ADO}_K(x,e_p)
  =
  \frac{
    e_p^t x^{(1-p)t}
  }{
    (1-x)(1+x^p)}
  \sum_{k=0}^{p-1}
  (-1)^k e_p^{\frac{2t+1}{2}k^2+\frac{2t-1}{2}k}
  x^{(2t+1)k}
  (1-e_p^{2k+1}x^2)
\end{equation*}
Inserting this into the definition of the CGP, we have
\begin{equation}\label{eq:CGPtorus}
  \begin{aligned}[t]
    \mathrm{CGP}(S^3(K),\lambda)
    & =
    \sum_{n=0}^{p-1}
    d^2(\lambda+2n) \, \mathrm{ADO}_K(e_n^{-(\lambda+2n+1)}, e_p)
    \\
    & =
    \frac{e_p^{(p-1) t \lambda}}{
      (e_{2p}^{p \lambda}-e_{2p}^{-p \lambda})^2
      (1+e_p^{-p \lambda})}
    \sum_{k=0}^{p-1}
    (-1)^k e_p^{\frac{2t+1}{2}k^2 + \frac{2t-1}{2}k
      -(\lambda+1)(2t+1)k}
    \\
    & \qquad \times
    \sum_{n=0}^{p-1}
    e_p^{-2n(k(2t+1)+t)}
    (e_p^{\lambda+2n+1}-1)
    (1-e_p^{2k-1-2\lambda-4n})
  \end{aligned}
\end{equation}
Observe that for $\lambda=0$ \eqref{eq:CGPtorus}
  coincides with~\eqref{eq:WRTtorus} if we forget about the normalizing factors in front of the sum.
Furthermore,
in both expressions 
the sum over $n$ can be computed using
~\eqref{sum_ep}. 
 The $a$th power of $e^{2n}_p$  has a nonzero contribution only if
$p\mid a$. Now for each $(e^{n}_p)^{kp}$ with
$k\in \Z$ in \eqref{eq:WRTtorus}, there will be a 
corresponding term
$(e^{\lambda +2n+1}_p)^{kp}$ contributing to
\eqref{eq:CGPtorus} with (up to normalization)
the same coefficient.
 Setting $(e^{k\lambda}_p)^p:=T^k$, we observe
 that (up to normalization) CGP is 
a Laurent polynomial in $T^{\pm 1}$, such that evaluated at $T=1$ it coincides
with WRT.

\appendix

\section{}
Let $K=\overline{T_{(2,5)}}$, a mirror image of
$T_{(2,5)}$.
The coefficient of the Habiro expansion  
\begin{equation}
  \label{a_for_25}
  a_n(q)
  =(-1)^n q^{\frac{1}{2} n^2 + \frac{3}{2} n -1}
  \sum_{k=1}^{n+1} q^{k^2}
  \begin{bmatrix}
    n+k \\
    2k-1
  \end{bmatrix}_q
\end{equation}
is obtained by setting $t=2$ in \eqref{Hab.coef}.
Using $k= \ell p + j$,
we get by~\eqref{fact:qbinom}
\begin{equation*}
  \begin{aligned}
    a_{ m p}(e_p)
    &=
    (-1)^{m} e_p^{-1}
    \left(
      e_p +
      \sum_{\ell=0}^{m-1}\sum_{j=1}^p
      e_p^{j^2}
      \begin{bmatrix}
        (m + \ell)p +j
        \\
        2(\ell p + j) -1
      \end{bmatrix}_{e_p}
    \right)
    \\
    & =
    (-1)^m
    \left(
      1+ \sum_{\ell=0}^{m-1}
      \begin{pmatrix}
        m+\ell \\
        2 \ell
      \end{pmatrix}
      +
      e_p^{-1}
      \sum_{\ell=0}^{m-1}
      \begin{pmatrix}
        m+\ell\\
        2\ell+1
      \end{pmatrix}
      \sum_{j=
        \lfloor \frac{p}{2} \rfloor + 1}^{p-1}
      e_p^{j^2}
      \begin{bmatrix}
        j
        \\
        2j-1-p
      \end{bmatrix}_{e_p}
    \right)
  \end{aligned}
\end{equation*}
Especially
we have
\begin{equation}
  a_p(e_p)
  =
  -2-
  \sum_{j=
    \lfloor \frac{p}{2} \rfloor +1
  }^{p-1}
  e_p^{j^2-1}
  \begin{bmatrix}
    j \\
    2j-1-p
  \end{bmatrix}_{e_p} .
\end{equation}
By definition~\eqref{a_for_25},  we have
\begin{gather*}
  a_{n-1}(1)
  =
  (-1)^{n-1} \sum_{\ell=0}^{n}
  \begin{pmatrix}
    n+\ell \\
    2\ell+ 1
  \end{pmatrix}
  \\
  a_{2m}(-1)
  =
  (-1)^m \left(
    1+\sum_{\ell=0}^{m-1}
    \begin{pmatrix}
      m+\ell
      \\
      2\ell
    \end{pmatrix}
  \right)
\end{gather*}
Combining these identities, we obtain
\begin{equation}
  \label{eq:25}
  a_{m p}(e_p)
  =
  a_{2 m} (-1) + a_{m-1}(1) \left( 2+ a_p(e_p) \right)
\end{equation}

\end{document}